\documentclass[10pt]{article}
\usepackage{amsmath,amsfonts,mathrsfs,amssymb}
\usepackage{indentfirst}

\begin{document}

\title{\textbf{Partial result of Yau's Conjecture of the first eigenvalue in unit sphere $\mathbb{S}^{n+1}(1)$}}
\author{\large{ Zhongyang Sun}\\
\small\footnotesize {}}
\date{}

\maketitle
\footnote[0]{2010 $Mathematics$ $Subject$ $Classification$. Primary 35P15; Secondary 53C40, 53C42.}
\footnote[0]{This work was supported by the National Natural Science Foundation of China (No. 11501500).}
\footnote[0]{{School of Mathematics Science, Huaibei Normal University, Huaibei, Anhui, 235000, China} }
\footnote[0]{{E-mail: sunzhongyang12@163.com} }

\vskip-.8cm {\normalsize $\mathbf{Abstract:}$ }{\small
In this paper, we partially solve Yau' Conjecture of the first eigenvalue of an embedded compact minimal hypersurface of unit sphere $\mathbb{S}^{n+1}(1)$, i.e., Corollary 1.2.
In particular, Corollary 1.3 proves that the condition $\int_{\Omega_{1}}|\nabla u|^{2}=(n+1)\int_{\Omega_{1}}u^{2}$ is naturally true and meaningful in Corollary 1.2.}

\vspace{6pt} {\normalsize $\mathbf{Keywords:}$ }{\small
Minimal hypersurface; Riemannian manifold ;
The first eigenvalue ; Unit sphere }

\section{Introduction and main results}

Let $M^{n}$ be an embedded compact orientable minimal hypersurface in an $(n+1)$-dimensional compact orientable Riemannian manifold $N$.

In 1982, Yau {\cite{[Y]}} proposed the following conjecture.

\vskip.2cm
{ {\noindent \bf Yau' Conjecture.}}
\textit{If $\lambda_{1}(M^{n})$ is the first eigenvalue of an embedded compact minimal hypersurface $M^{n}$ of unit sphere $\mathbb{S}^{n+1}(1)$,
the standard $(n+1)$-sphere of sectional curvature 1, then $\lambda_{1}(M^{n})=n$.}

In 1983, Choi-Wang \cite{[CW]} showed that if $M^{n}$ is an embedded compact minimal hypersurface of unit sphere $\mathbb{S}^{n+1}(1)$,
then $\lambda_{1}(M^{n})\geqslant \frac{n}{2}$, where $\lambda_{1}(M^{n})$ is the first eigenvalue of the Laplacian of $M^{n}$.

In the paper, in order to prove Yau' Conjecture, first we obtain the following result in Riemannian manifold $N$.
\vskip.2cm
{\noindent \bf Theorem 1.1.}
\textit{Let $M^{n}$ be an embedded compact orientable minimal hypersurface in an $(n+1)$-dimensional compact orientable Riemannian manifold $N$.
Suppose that $\lambda_{1}(M^{n})$ is the first eigenvalue of the Laplacian of $M^{n}$.
If the Ricci curvature of $N$ is bounded below by a positive constant $k$, then
$$
\frac{k}{2}+\frac{n}{2}\delta_{1}(u)\leqslant\lambda_{1}(M^{n})\leqslant\frac{k}{2}+\frac{n}{2}\delta_{2}(u),
$$
where
$$
0<\delta_{1}(u)=\frac{\int_{\Omega_{1}}|\nabla u|^{2}
\left(1-\sqrt{1-\frac{(n+1)k}{n}\frac{\int_{\Omega_{1}}u^{2}}{\int_{\Omega_{1}}|\nabla u|^{2}}}\right)^{2}}
{(n+1)\int_{\Omega_{1}}u^{2}}\leqslant\frac{k}{n},
$$
$$
\delta_{2}(u)=\frac{\int_{\Omega_{1}}|\nabla u|^{2}
\left(1+\sqrt{1-\frac{(n+1)k}{n}\frac{\int_{\Omega_{1}}u^{2}}{\int_{\Omega_{1}}|\nabla u|^{2}}}\right)^{2}}
{(n+1)\int_{\Omega_{1}}u^{2}}\geqslant\frac{k}{n},
$$
and $u$ is the solution of the Dirichlet problem such that
$$
\left\{
  \begin{array}{ll}
    \Delta u=0,~~~~~~~~in~\Omega_{1};\\
    u=f,~~~~~~~~~~in~\partial\Omega_{1}=M^{n},
  \end{array}
\right.
$$
which is called $\mathrm{Dirichlet}$ $\mathrm{problem}$ $\mathrm{(1)}$.}

\vskip.2cm
We see that $0<\delta_{1}(u)\leqslant\frac{k}{n}$, $\delta_{2}(u)\geqslant\frac{k}{n}$ in Theorem 1.1.
In particular, we have $\delta_{1}(u)=\delta_{2}(u)=\frac{k}{n}$ iff
$\int_{\Omega_{1}}|\nabla u|^{2}=\frac{(n+1)k}{n}\int_{\Omega_{1}}u^{2}$.
When we take $\int_{\Omega_{1}}|\nabla u|^{2}=\frac{(n+1)k}{n}\int_{\Omega_{1}}u^{2}$ in Theorem 1.1, we have $\lambda_{1}(M^{n})=k$.
Note that the Ricci curvature of unit sphere $\mathbb{S}^{n+1}(1)$ is $n$ and an embedded compact hypersurface of unit sphere $\mathbb{S}^{n+1}(1)$ is already orientable.
We obtain Corollary 1.2 which partially solves Yau' Conjecture.

\vskip.2cm
{\noindent \bf Corollary 1.2.}
\textit{Let $M^{n}$ be an embedded compact minimal hypersurface of unit sphere $\mathbb{S}^{n+1}(1)$.
Suppose that $\lambda_{1}(M^{n})$ is the first eigenvalue of the Laplacian of $M^{n}$.
If the solution $u$ of $\mathrm{Dirichlet}$ $\mathrm{problem}$ $\mathrm{(1)}$ satisfies $\int_{\Omega_{1}}|\nabla u|^{2}=(n+1)\int_{\Omega_{1}}u^{2}$, then
$$
\lambda_{1}(M^{n})=n.
$$}

According to (3.3) in the proof of Theorem 1.1, we know that $\int_{\Omega_{1}}|\nabla u|^{2}\geqslant\frac{(n+1)k}{n}\int_{\Omega_{1}}u^{2}$ is naturally true
and meaningful in Theorem 1.1.
Note that the Ricci curvature of unit sphere $\mathbb{S}^{n+1}(1)$ is $n$. Hence, we get the following corollary.

\vskip.2cm
{\noindent \bf Corollary 1.3.}
\textit{Let $M^{n}$ be an embedded compact minimal hypersurface of unit sphere $\mathbb{S}^{n+1}(1)$.
Suppose that $\lambda_{1}(M^{n})$ is the first eigenvalue of the Laplacian of $M^{n}$.
Then
$$
\int_{\Omega_{1}}|\nabla u|^{2}\geqslant(n+1)\int_{\Omega_{1}}u^{2},
$$
where $u$ is the solution of $\mathrm{Dirichlet}$ $\mathrm{problem}$ $\mathrm{(1)}$.}

\vskip.2cm
{\noindent \bf Remark 1.4.}
Corollary 1.3 proves that the condition $\int_{\Omega_{1}}|\nabla u|^{2}=(n+1)\int_{\Omega_{1}}u^{2}$ is naturally true and meaningful in Corollary 1.2.

\vskip.2cm
From Lemma 2.1 and $\int_{\partial\Omega_{1}}h(\overline{\nabla} u,\overline{\nabla} u)\geqslant0$ in the proof of Theorem 1.1,
we can conclude the following corollary.

\vskip.2cm
{\noindent \bf Corollary 1.5.}
\textit{Let $M^{n}$ be an embedded compact orientable minimal hypersurface in an $(n+1)$-dimensional compact orientable Riemannian manifold $N$.
Suppose that $\lambda_{1}(M^{n})$ is the first eigenvalue of the Laplacian of $M^{n}$.
If the Ricci curvature of $N$ is bounded below by a positive constant $k$, then
$$
\lambda_{1}(M^{n})>\frac{k}{2}.
$$}

Note that the Ricci curvature of unit sphere $\mathbb{S}^{n+1}(1)$ is $n$.
Hence, from Corollary 1.5 we obtain Corollary 1.6 which improves the result $\lambda_{1}(M^{n})\geqslant\frac{n}{2}$ of \cite{[CW]}.

\vskip.2cm
{\noindent \bf Corollary 1.6.}
\textit{Let $M^{n}$ be an embedded compact minimal hypersurface of unit sphere $\mathbb{S}^{n+1}(1)$.
Suppose that $\lambda_{1}(M^{n})$ is the first eigenvalue of the Laplacian of $M^{n}$.
Then
$$
\lambda_{1}(M^{n})>\frac{n}{2}.
$$}

Finally, by Corollaries 1.2 and 1.3, what we want to know is whether the following problem is true.

\vskip.2cm
{\noindent \bf Problem 1.7.}
\textit{Let $M^{n}$ be an embedded compact minimal hypersurface of unit sphere $\mathbb{S}^{n+1}(1)$.
Suppose that $\lambda_{1}(M^{n})$ is the first eigenvalue of the Laplacian of $M^{n}$.
If the solution $u$ of $\mathrm{Dirichlet}$ $\mathrm{problem}$ $\mathrm{(1)}$ satisfies $\int_{\Omega_{1}}|\nabla u|^{2}>(n+1)\int_{\Omega_{1}}u^{2}$,
then $\lambda_{1}(M^{n})=n$~$\mathrm{?}$}

\vskip.2cm
{\noindent \bf Remark 1.8.}
If Problem 1.7 is true, combining Corollary 1.2, Corollary 1.3 and Problem 1.7, then Yau' Conjecture is true.

\section{Preliminaries}

Let $\Omega_{1}$ is a Riemannian manifold of dimensional $(n+1)$ with smooth boundary $\partial\Omega_{1}=M^{n}$.
Let $u$ be a function defined on $\Omega_{1}$ which is smooth up to $\partial\Omega_{1}$.
The symbols $\Delta u$ and $\nabla u$ will be respectively the Laplacian and the gradient of $u$
with respect to the induced Riemannian metric on $\Omega_{1}$
while $\overline{\Delta} u$ and $\overline{\nabla} u$ will be the Laplacian
the gradient of $u$ (defined on $\partial\Omega_{1}$) with respect to the induced Riemannian metric on $\partial\Omega_{1}$.
For $x\in \Omega_{1}$ and $X,Y\in T_{x}\Omega_{1}$,
we can define the Hessian tensor $(D^{2}u)(X,Y)=X(Yu)-(\nabla_{X}Y)u$,
where $\nabla_{X}Y$ is the coviariant derivative of the Riemannian connection of $\Omega_{1}$.
The covariant derivative of the Riemannian connection of $\partial\Omega_{1}=M^{n}$ is given by $\overline{\nabla}_{X}Y$.

Suppose that $\{e_{1},e_{2},\cdots,e_{n},e_{n+1}\}$ is a local orthonormal frame such that
at $x\in \partial\Omega_{1}$, $e_{1},e_{2}$,$\cdots,e_{n},$ are tangent to $\partial\Omega_{1}$
and $e_{n+1}$ is the outward normal vector.
Let $h$ be the second fundamental form, $h(\upsilon,\omega)=<\nabla_{\upsilon}e_{n},\omega>$,
where $\upsilon,\omega$ are vectors tangent to $\partial\Omega_{1}=M^{n}$, and $H$ be the mean curvature,
i.e., $H=\sum^{n}_{i=1}\frac{h(e_{i},e_{i})}{n}$.

\vskip.2cm
We need the following lemmas will play a crucial role in the proof of Theorem 1.1.

\vskip.2cm
{\noindent \bf Lemma 2.1.}
\textit{Let $M^{n}$ be an embedded compact orientable minimal hypersurface in an $(n+1)$-dimensional compact orientable Riemannian
manifold $N$.
Suppose that $\lambda_{1}(M^{n})$ is the first eigenvalue of the Laplacian of $M^{n}$.
If the Ricci curvature of $N$ is bounded below by a positive constant $k$, then
$$
\big(2\lambda_{1}(M^{n})-k\big)\int_{\Omega_{1}}|\nabla u|^{2}-
\int_{\partial\Omega_{1}}h(\overline{\nabla} u,\overline{\nabla} u)>0,
$$
where $u$ is the solution of $\mathrm{Dirichlet}$ $\mathrm{problem}$ $\mathrm{(1)}$.}

\vskip.2cm
{\noindent \bf Proof.}
Since the Ricci curvature of $N$ is strictly positive, the first Betti number of $N$ must be zero.
Since $M^{n}$ and $N$ are orientable, by looking at the exact sequences of homology groups,
we can see that $M^{n}$ divides $N$ into components $\Omega_{1}$ and $\Omega_{2}$
such that $\partial\Omega_{1}=\partial\Omega_{2}=M^{n}$.

Let $f$ be the first eigenfunction of $M^{n}$, i.e.,
$$
\overline{\Delta}f+\lambda_{1}(M^{n})f=0.
\eqno(2.1)
$$
Let $u$ be the solution of the Dirichlet problem such that
$$
\left\{
  \begin{array}{ll}
    \Delta u=0,~~~~~~~~\textrm{in}~~\Omega_{1};\\
    u=f,~~~~~~~~~~\textrm{in}~~\partial\Omega_{1}=M^{n},
  \end{array}
\right.
\eqno(2.2)
$$
which is called Dirichlet problem (1).

So $u$ is a function defined on $\Omega_{1}$ smooth up to $\partial\Omega_{1}$.
Then we have
$$
\Delta u=\sum_{i=1}^{n+1}D^{2}u(e_{i},e_{i})=\sum_{i=1}^{n+1}u_{ii},
\eqno(2.3)
$$
where $u_{ij}=D^{2}u(e_{i},e_{j})$, ~~~~~~~~~$i,j=1,\cdots,n+1$.

When $i\neq n+1$ and $x\in \partial\Omega_{1}$, we get
$$
\nabla_{e_{i}}e_{i}=\overline{\nabla}_{e_{i}}e_{i}-h_{ii}e_{n+1},
\eqno(2.4)
$$
where $h_{ij}=h(e_{i},e_{j})$.

Hence when $x\in \partial\Omega_{1}=M^{n}$, we can conclude from (2.2), (2.3) and (2.4) that
$$
\Delta u=u_{(n+1)(n+1)}+\overline{\Delta}f+\sum_{i=1}^{n}h_{ii}e_{n+1}(u)
=u_{(n+1)(n+1)}+\overline{\Delta}f+nHu_{n+1},
\eqno(2.5)
$$
where $u_{n+1}=<\nabla u, e_{n+1}>$ and $H$ is the mean curvature of $M^{n}$.

For $x\in\partial\Omega_{1}$, we know from (2.1) and (2.5) that
$$
u_{(n+1)(n+1)}=\lambda_{1}(M^{n})f-nHu_{n+1}.
\eqno(2.6)
$$
For $x\in\Omega_{1}$, note that the fact
$\Delta|\nabla u|^{2}
=2\sum_{i,j=1}^{n+1}u_{ij}^{2}+2\sum_{i,j=1}^{n+1}R_{ij}u_{i}u_{j}+2\sum_{i=1}^{n+1}u_{i}(\Delta u)_{i}$
can be found in \cite{[SY]}. Since $\Delta u=0$, we have
$$
\Delta|\nabla u|^{2}=2\sum_{i,j=1}^{n+1}u_{ij}^{2}+2\sum_{i,j=1}^{n+1}R_{ij}u_{i}u_{j},
\eqno(2.7)
$$
where $R_{ij}=\mathrm{Ric}(e_{i},e_{j})$.

Since Ricci curvature of $N$ is bounded by $k$, by (2.7) we have
$$
\Delta|\nabla u|^{2}\geqslant2|D^{2}u|^{2}+2k|\nabla u|^{2}.
\eqno(2.8)
$$
Then integrating (2.8) we have
$$
\int_{\Omega_{1}}\Delta|\nabla u|^{2}\geqslant2\int_{\Omega_{1}}|D^{2}u|^{2}+2k\int_{\Omega_{1}}|\nabla u|^{2}.
\eqno(2.9)
$$
When $i\neq n+1$, we have
$$
u_{i(n+1)}=D^{2}u(e_{i},e_{n+1})=e_{i}(e_{n+1}u)-(\nabla_{e_{i}}e_{n+1})u=e_{i}(u_{n+1})-\sum_{j=1}^{n}h_{ij}u_{j}.
\eqno(2.10)
$$
By using the Stokes theorem, we can conclude from (2.1),(2.6) and (2.10) that
$$
\begin{aligned}
\int_{\Omega_{1}}\Delta|\nabla u|^{2}
=&2\int_{\partial\Omega_{1}}\sum_{i=1}^{n}u_{i}u_{i(n+1)}+2\int_{\partial\Omega_{1}}u_{n+1}u_{(n+1)(n+1)}\\
=&2\int_{\partial\Omega_{1}}\overline{\nabla}f\cdot\overline{\nabla}u_{n+1}-2\int_{\partial\Omega_{1}}\sum_{i,j=1}^{n}h_{ij}u_{i}u_{j}+
2\int_{\partial\Omega_{1}}u_{n+1}u_{(n+1)(n+1)}\\
=&-2\int_{\partial\Omega_{1}}u_{n+1}\overline{\Delta}f-2\int_{\partial\Omega_{1}}h(\overline{\nabla} u,\overline{\nabla} u)
+2\int_{\partial\Omega_{1}}u_{n+1}u_{(n+1)(n+1)}\\
=&4\lambda_{1}(M^{n})\int_{\partial\Omega_{1}}u_{n+1}f-2\int_{\partial\Omega_{1}}h(\overline{\nabla} u,\overline{\nabla} u)
-2n\int_{\partial\Omega_{1}}Hu_{n+1}^{2}.
\end{aligned}
\eqno(2.11)
$$
Combining (2.2) and Stokes theorem, we have
$$
\int_{\Omega_{1}}|\nabla u|^{2}=-\int_{\Omega_{1}}u\Delta u+\int_{\partial\Omega_{1}}uu_{n+1}
=\int_{\partial\Omega_{1}}u_{n+1}f.
\eqno(2.12)
$$
From (2.11) and (2.12), we have
$$
\int_{\Omega_{1}}\Delta|\nabla u|^{2}=4\lambda_{1}(M^{n})\int_{\Omega_{1}}|\nabla u|^{2}
-2\int_{\partial\Omega_{1}}h(\overline{\nabla} u,\overline{\nabla} u)-2n\int_{\partial\Omega_{1}}Hu_{n+1}^{2}.
\eqno(2.13)
$$
Since $M^{n}$ is minimal, by (2.9) and (2.13), we have
$$
\big(2\lambda_{1}(M^{n})-k\big)\int_{\Omega_{1}}|\nabla u|^{2}\geqslant\int_{\partial\Omega_{1}}h(\overline{\nabla} u,\overline{\nabla} u)
+\int_{\Omega_{1}}|D^{2}u|^{2}.
\eqno(2.14)
$$
We claim that
$$
\int_{\Omega_{1}}|D^{2}u|^{2}\neq0.
\eqno(2.15)
$$
Otherwise, for all $1\leqslant i,j\leqslant n+1$, we have $u_{ij}=0$ on $\Omega_{1}$.
Since $u$ is smooth up to $\partial\Omega_{1}$, for all $1\leqslant i,j\leqslant n$, we have $f_{ij}=0$ on $M^{n}$,
which implies that $\overline{\Delta}f=0$ which is impossible since $f$ is the first eigenfunction of $M^{n}$.
Thus, our claim is true. Hence from (2.15) we get
$$
\int_{\Omega_{1}}|D^{2}u|^{2}>0.
\eqno(2.16)
$$
From (2.14) and (2.16), we have
$$
\big(2\lambda_{1}(M^{n})-k\big)\int_{\Omega_{1}}|\nabla u|^{2}-\int_{\partial\Omega_{1}}h(\overline{\nabla} f,\overline{\nabla} f)>0.
\eqno(2.17)
$$
This completes the proof of Lemma 2.1.
\hfill$\Box$

\vskip.2cm
{\noindent \bf Lemma 2.2.}
\textit{Let $M^{n}$ be an embedded compact orientable minimal hypersurface in an $(n+1)$-dimensional compact orientable Riemannian manifold $N$.
Suppose that $\lambda_{1}(M^{n})$ is the first eigenvalue of the Laplacian of $M^{n}$.
If the Ricci curvature of $N$ is bounded below by a positive constant $k$, then
$$
\begin{aligned}
Q(t)=&\left[\big(2\lambda_{1}(M^{n})-k\big)\int_{\Omega_{1}}|\nabla u|^{2}-\int_{\partial\Omega_{1}}h(\overline{\nabla}f,\overline{\nabla}f)\right]\cdot t^{2}
+2\lambda_{1}(M^{n})\int_{\Omega_{1}}u^{2}\cdot t\\
&+\frac{n}{n+1}\int_{\Omega_{1}}u^{2}\geqslant0,~~~~~~~\forall~t\in \mathbb{R},
\end{aligned}
\eqno(2.18)
$$
where $u$ is the solution of $\mathrm{Dirichlet}$ $\mathrm{problem}$ $\mathrm{(1)}$.}

\vskip.2cm
{\noindent \bf Proof.}
Let $z=\frac{\partial u}{\partial \nu}=u_{\nu}$ be the normal outward derivative of $u$.
We know that the Reilly formula:
$$
\int_{\Omega_{1}}(\Delta u)^{2}=\int_{\Omega_{1}}|D^{2}u|^{2}+\int_{\Omega_{1}}\mathrm{Ric}(\nabla u,\nabla u)
+\int_{\partial\Omega_{1}}2z\overline{\Delta}u+\int_{\partial\Omega_{1}}h(\overline{\nabla}u,\overline{\nabla}u)
+\int_{\partial\Omega_{1}}nHz^{2}.
\eqno(2.19)
$$
Case 1: $t=0$, then (2.18) is naturally true.
Case 2: $t\neq0$, we will consider the following Dirichlet problem
$$
\left\{
  \begin{array}{ll}
    \Delta g=u,~~~~~~~~\textrm{in}~~\Omega_{1};\\
    g=tf,~~~~~~~~~\textrm{in}~~\partial\Omega_{1}=M^{n}.
  \end{array}
\right.
\eqno(2.20)
$$
Combining Green formula and (2.20), we have
$$
\left\{
  \begin{array}{ll}
    \int_{\partial\Omega_{1}}f\frac{\partial u}{\partial \nu}=\int_{\Omega_{1}}|\nabla u|^{2};\\
    t\int_{\partial\Omega_{1}}f\frac{\partial u}{\partial \nu}=\int_{\Omega_{1}}<\nabla u,\nabla g>;\\
    \int_{\partial\Omega_{1}}f\frac{\partial g}{\partial \nu}=\int_{\Omega_{1}}u^{2}+\int_{\Omega_{1}}<\nabla u,\nabla g>.
  \end{array}
\right.
\eqno(2.21)
$$
From (2.21), we have
$$
\int_{\Omega_{1}}<\nabla u,\nabla g>=t\int_{\Omega_{1}}|\nabla u|^{2}.
\eqno(2.22)
$$
From (2.22) and Cauchy-Schwarz inequality, we have
$$
\int_{\Omega_{1}}|\nabla g|^{2}\geqslant t^{2}\int_{\Omega_{1}}|\nabla u|^{2}.
\eqno(2.23)
$$
From the third equation in (2.21) and (2.22), we get
$$
t\int_{\partial\Omega_{1}}f\frac{\partial g}{\partial \nu}=t\int_{\Omega_{1}}u^{2}+t^{2}\int_{\Omega_{1}}|\nabla u|^{2}.
\eqno(2.24)
$$
Since $M^{n}$ is minimal and $|D^{2}g|^{2}\geqslant \frac{1}{n+1}(\Delta g)^{2}$, applying (2.19) to $g$, we obtain
$$
\frac{n}{n+1}\int_{\Omega_{1}}(\Delta g)^{2}\geqslant k\int_{\Omega_{1}}|\nabla g|^{2}
+2\int_{\partial\Omega_{1}}\frac{\partial g}{\partial\nu}\overline{\Delta}(tf)+\int_{\partial\Omega_{1}}h(\overline{\nabla}g,\overline{\nabla}g).
\eqno(2.25)
$$
On the other hand, combining (2.1), (2.23), (2.24), (2.25) and $\Delta g=u$, we have
$$
\frac{n}{n+1}\int_{\Omega_{1}}u^{2}\geqslant kt^{2}\int_{\Omega_{1}}|\nabla u|^{2}
-2\lambda_{1}(M^{n})\left[t\int_{\Omega_{1}}u^{2}+t^{2}\int_{\Omega_{1}}|\nabla u|^{2}\right]
+t^{2}\int_{\partial\Omega_{1}}h(\overline{\nabla}f,\overline{\nabla}f).
\eqno(2.26)
$$
Hence, we have
$$
\begin{aligned}
Q(t)=&\left[\big(2\lambda_{1}(M^{n})-k\big)\int_{\Omega_{1}}|\nabla u|^{2}-\int_{\partial\Omega_{1}}h(\overline{\nabla}f,\overline{\nabla}f)\right]\cdot t^{2}
+2\lambda_{1}(M^{n})\int_{\Omega_{1}}u^{2}\cdot t\\
&+\frac{n}{n+1}\int_{\Omega_{1}}u^{2}\geqslant0,~~~~~~~\forall~t\in\mathbb{R}.
\end{aligned}
$$
This completes the proof of Lemma 2.2.
\hfill$\Box$

\section{Proof of Theorem 1.1}

\vskip.2cm
{\noindent \bf Proof of Theorem 1.1.}
Combining (2.17) and (2.18), we know that the discriminant of $Q(t)$ is non-positive, i.e.,
$$
\frac{n+1}{n}\lambda_{1}^{2}(M^{n})\frac{\int_{\Omega_{1}}u^{2}}{\int_{\Omega_{1}}|\nabla u|^{2}}
-2\lambda_{1}(M^{n})+k+\frac{\int_{\partial\Omega_{1}}h(\overline{\nabla} f,\overline{\nabla} f)}{\int_{\Omega_{1}}|\nabla u|^{2}}\leqslant0.
\eqno(3.1)
$$
Since $\int_{\partial\Omega_{1}}h(\overline{\nabla} u,\overline{\nabla} u)=\int_{M^{n}}h(\overline{\nabla} f,\overline{\nabla} f)$ and
the outward normal vector of $\partial\Omega_{2}$ is $-e_{n}$,
we have
$$
\int_{\partial\Omega_{2}}h(\overline{\nabla} u,\overline{\nabla} u)=-\int_{\partial\Omega_{1}}h(\overline{\nabla} u,\overline{\nabla} u).
$$
Hence we can assume that $\int_{\partial\Omega_{1}}h(\overline{\nabla} u,\overline{\nabla} u)\geqslant0$; otherwise,
we work with $\Omega_{2}$ rather with $\Omega_{1}$.
Since $\int_{\partial\Omega_{1}}h(\overline{\nabla}f,\overline{\nabla}f)\geqslant0$, from (3.1) we have
$$
\frac{n+1}{n}\lambda_{1}^{2}(M^{n})\frac{\int_{\Omega_{1}}u^{2}}{\int_{\Omega_{1}}|\nabla u|^{2}}
-2\lambda_{1}(M^{n})+k\leqslant0.
\eqno(3.2)
$$
Combining (2.17), $\int_{\partial\Omega_{1}}h(\overline{\nabla}f,\overline{\nabla}f)\geqslant0$ and $\int_{\Omega_{1}}|\nabla u|^{2}>0$,
we have $\lambda_{1}(M^{n})>\frac{k}{2}$.
When we take $t=-\frac{n}{(n+1)k}$ in (2.18),
combining $\lambda_{1}(M^{n})>\frac{k}{2}$, $\int_{\partial\Omega_{1}}h(\overline{\nabla}f,\overline{\nabla}f)\geqslant0$ and $\int_{\Omega_{1}}|\nabla u|^{2}>0$,
we get
$$
1-\frac{(n+1)k}{n}\frac{\int_{\Omega_{1}}u^{2}}{\int_{\Omega_{1}}|\nabla u|^{2}}\geqslant0.
\eqno(3.3)
$$
From (3.2) and (3.3), we have
$$
\begin{aligned}
\lambda_{1}(M^{n})\geqslant&\frac{n\int_{\Omega_{1}}|\nabla u|^{2}-n\int_{\Omega_{1}}|\nabla u|^{2}\sqrt{1-\frac{(n+1)k}{n}\frac{\int_{\Omega_{1}}u^{2}}
{\int_{\Omega_{1}}|\nabla u|^{2}}}}{(n+1)\int_{\Omega_{1}}u^{2}}\\
=&\frac{k(n+1)\int_{\Omega_{1}}u^{2}}{2n(n+1)\int_{\Omega_{1}}u^{2}}\\
&+\frac{2n\int_{\Omega_{1}}|\nabla u|^{2}-k(n+1)\int_{\Omega_{1}}u^{2}
-2n\int_{\Omega_{1}}|\nabla u|^{2}\sqrt{1-\frac{(n+1)k}{n}\frac{\int_{\Omega_{1}}u^{2}}
{\int_{\Omega_{1}}|\nabla u|^{2}}}}{2(n+1)\int_{\Omega_{1}}u^{2}}\\
=&\frac{k}{2}+\frac{n}{2}\delta_{1}(u)
\end{aligned}
\eqno(3.4)
$$
and
$$
\begin{aligned}
\lambda_{1}(M^{n})\leqslant&\frac{n\int_{\Omega_{1}}|\nabla u|^{2}+n\int_{\Omega_{1}}|\nabla u|^{2}\sqrt{1-\frac{(n+1)k}{n}\frac{\int_{\Omega_{1}}u^{2}}
{\int_{\Omega_{1}}|\nabla u|^{2}}}}{(n+1)\int_{\Omega_{1}}u^{2}}\\
=&\frac{k(n+1)\int_{\Omega_{1}}u^{2}}{2n(n+1)\int_{\Omega_{1}}u^{2}}\\
&+\frac{2n\int_{\Omega_{1}}|\nabla u|^{2}-k(n+1)\int_{\Omega_{1}}u^{2}
+2n\int_{\Omega_{1}}|\nabla u|^{2}\sqrt{1-\frac{(n+1)k}{n}\frac{\int_{\Omega_{1}}u^{2}}
{\int_{\Omega_{1}}|\nabla u|^{2}}}}{2(n+1)\int_{\Omega_{1}}u^{2}}\\
=&\frac{k}{2}+\frac{n}{2}\delta_{2}(u),
\end{aligned}
\eqno(3.5)
$$
where
$$
\delta_{1}(u)=\frac{\int_{\Omega_{1}}|\nabla u|^{2}
\left(1-\sqrt{1-\frac{(n+1)k}{n}\frac{\int_{\Omega_{1}}u^{2}}{\int_{\Omega_{1}}|\nabla u|^{2}}}\right)^{2}}
{(n+1)\int_{\Omega_{1}}u^{2}}
\eqno(3.6)
$$
and
$$
\delta_{2}(u)=\frac{\int_{\Omega_{1}}|\nabla u|^{2}
\left(1+\sqrt{1-\frac{(n+1)k}{n}\frac{\int_{\Omega_{1}}u^{2}}{\int_{\Omega_{1}}|\nabla u|^{2}}}\right)^{2}}
{(n+1)\int_{\Omega_{1}}u^{2}}.
\eqno(3.7)
$$
Finally we will prove $0<\delta_{1}(u)\leqslant\frac{k}{n}$, $\delta_{2}(u)\geqslant\frac{k}{n}$.
Firstly, from (3.6) we have $\delta_{1}(u)\geqslant0$. But $\delta_{1}(u)\neq0$, otherwise, from (3.6) we can conclude that
$1-\sqrt{1-\frac{(n+1)k}{n}\frac{\int_{\Omega_{1}}u^{2}}{\int_{\Omega_{1}}|\nabla u|^{2}}}=0$,
which implies that $\int_{\Omega_{1}}u^{2}=0$ which is impossible.
Hence we get $\delta_{1}(u)>0$ and $1-\frac{(n+1)k}{n}\frac{\int_{\Omega_{1}}u^{2}}{\int_{\Omega_{1}}|\nabla u|^{2}}\neq1$.
Combining (3.3), we know that
$$
0\leqslant1-\frac{(n+1)k}{n}\frac{\int_{\Omega_{1}}u^{2}}{\int_{\Omega_{1}}|\nabla u|^{2}}<1.
\eqno(3.8)
$$
Obviously, $\delta_{1}(u)$ can be written as
$$
\delta_{1}(u)
=\frac{2\int_{\Omega_{1}}|\nabla u|^{2}
\left(1-\frac{(n+1)k}{n}\frac{\int_{\Omega_{1}}u^{2}}{\int_{\Omega_{1}}|\nabla u|^{2}}-\sqrt{1-\frac{(n+1)k}{n}\frac{\int_{\Omega_{1}}u^{2}}
{\int_{\Omega_{1}}|\nabla u|^{2}}}\right)+\frac{(n+1)k}{n}\int_{\Omega_{1}}u^{2}}
{(n+1)\int_{\Omega_{1}}u^{2}}.
\eqno(3.9)
$$
From (3.8), we obtain
$$
1-\frac{(n+1)k}{n}\frac{\int_{\Omega_{1}}u^{2}}{\int_{\Omega_{1}}|\nabla u|^{2}}
-\sqrt{1-\frac{(n+1)k}{n}\frac{\int_{\Omega_{1}}u^{2}}{\int_{\Omega_{1}}|\nabla u|^{2}}}\leqslant0.
\eqno(3.10)
$$
We can conclude from (3.9) and (3.10) that
$$
\delta_{1}(u)\leqslant\frac{k}{n}.
$$
Secondly, from (3.7) and (3.8) we obtain
$$
\begin{aligned}
\delta_{2}(u)=&\frac{\int_{\Omega_{1}}|\nabla u|^{2}
\left(1+\sqrt{1-\frac{(n+1)k}{n}\frac{\int_{\Omega_{1}}u^{2}}{\int_{\Omega_{1}}|\nabla u|^{2}}}\right)^{2}}
{(n+1)\int_{\Omega_{1}}u^{2}}\\
=&\frac{2\int_{\Omega_{1}}|\nabla u|^{2}
\left(1-\frac{(n+1)k}{n}\frac{\int_{\Omega_{1}}u^{2}}{\int_{\Omega_{1}}|\nabla u|^{2}}+\sqrt{1-\frac{(n+1)k}{n}\frac{\int_{\Omega_{1}}u^{2}}
{\int_{\Omega_{1}}|\nabla u|^{2}}}\right)+\frac{(n+1)k}{n}\int_{\Omega_{1}}u^{2}}
{(n+1)\int_{\Omega_{1}}u^{2}}\\
\geqslant&\frac{k}{n}.
\end{aligned}
$$
So the proof of Theorem 1.1 is finished.

\hfill$\Box$

\section{Applications}

In 1980, Yang and Yau \cite{[YY]} obtained the following result.
\vskip.2cm
{\noindent \bf Theorem 4.1 (\cite{[YY]}).}
\textit{If $M^{2}$ is an orientable Riemannian surface of genus $\emph{g}$ with area $A$,
then
$$
\lambda_{1}(M^{2})\leqslant 8\pi(\emph{g}+1)A^{-1}.
$$}

From Corollary 1.5 and Theorem 4.1, we obtain the following theorem.

\vskip.2cm
{\noindent \bf Theorem 4.2.}
\textit{Let $M^{2}$ be an embedded compact orientable minimal surface of genus $\emph{g}$ with area $A$ in a $3$-dimensional compact orientable Riemannian manifold $N$. If the Ricci curvature of $N$ is bounded below by a positive constant $k$, then
$$
A<\frac{16\pi(\emph{g}+1)}{k}.
$$}

Note that the Ricci curvature of unit sphere $\mathbb{S}^{n+1}(1)$ is $n$. From Theorem 4.2,  we get the following corollary.

\vskip.2cm
{\noindent \bf Corollary 4.3.}
\textit{Let $M^{2}$ be an embedded compact minimal surface of genus $\emph{g}$ with area $A$ in unit sphere $\mathbb{S}^{3}(1)$. Then
$$
A<8\pi(\emph{g}+1).
$$}

From Corollary 1.2 and Theorem 4.1, we get the following theorem.

\vskip.2cm
{\noindent \bf Theorem 4.4.}
\textit{Let $M^{2}$ be an embedded compact minimal surface of genus $\emph{g}$ with area $A$ in unit sphere $\mathbb{S}^{3}(1)$.
If the solution $u$ of $\mathrm{Dirichlet}$ $\mathrm{problem}$ $\mathrm{(1)}$ satisfies $\int_{\Omega_{1}}|\nabla u|^{2}=3\int_{\Omega_{1}}u^{2}$, then
$$
A\leqslant 4\pi(\emph{g}+1).
$$}

Naturally, we propose the following conjecture.

\vskip.2cm
{\noindent \bf Conjecture 4.5.}
\textit{Let $M^{2}$ be an embedded compact minimal surface of genus $\emph{g}$ with area $A$ in unit sphere $\mathbb{S}^{3}(1)$.
Then
$$
A\leqslant 4\pi(\emph{g}+1).
$$}

\bibliographystyle{amsplain}

\end{document}